\documentclass[12pt,reqno]{article}
\usepackage{amssymb,amscd,amsmath,amsthm,color}
\newtheorem{theorem}{Theorem}[section]
\newtheorem{lemma}[theorem]{Lemma}
\newtheorem{cor}[theorem]{Corollary}
\newtheorem{prop}[theorem]{Proposition}
\usepackage{enumerate,amssymb}
\theoremstyle{definition}

\theoremstyle{remark}
\newtheorem{remark}[theorem]{Remark}

\numberwithin{equation}{section}

\def\bM{\mathbb{M}}

\begin{document}
\baselineskip=15pt

\title{ Unitary orbits  of Hermitian operators with convex or concave functions }

\author{Jean-Christophe Bourin and Eun-Young Lee} 


\date{ }


\maketitle

\begin{abstract}

\noindent
This short but self-contained survey  presents a number of elegant   matrix/operator  inequalities for general convex or concave functions,  obtained with a unitary orbit technique.   Jensen,  sub or super-additivity  type inequalities   are considered. Some of them are substitutes to classical  inequalities (Choi, Davis, Hansen-Pedersen) for operator convex or concave functions. Various trace, norm and determinantal inequalities are derived. Combined with an interesting  decomposition for positive semi-definite matrices, several results for partitioned matrices are also obtained.
\end{abstract}
\maketitle

{\small\noindent
 
AMS subjects classification 2000: Primary 15A60, 47A30, 47A60}

\vskip 15pt

\section{Introduction}

\noindent
The functional analytic aspect of 
 Matrix Analysis is evident when matrices or operators are considered as non-commutative numbers, sequences or functions.  In particular, a significant part of this theory consists in  establishing theorems for Hermitian matrices regarded
 as generalized real numbers or functions.  Two classical trace inequalities may illustrate quite well this assertion. Given two Hermitian matrices $A$, $B$ and a concave function $f(t)$ defined on the real line,
\begin{equation}\label{VN}
{\mathrm{Tr\,}} f\left(\frac{ A+B}{2}\right) \ge {\mathrm{Tr\,}}\frac{f(A) +f(B)}{2}
\end{equation}
and, if further $f(0)\ge 0$ and both $A$ and $B$ are positive semi-definite,
\begin{equation}\label{Rot}
{\mathrm{Tr\,}} f(A+B) \le {\mathrm{Tr\,}} f(A)+f(B)
\end{equation}
The first inequality goes back to von-Neumann in the 1920's, the second is more subtle and has been proved only in 1969 by Rotfel'd \cite{Rot}. These trace inequalities are matrix versions of  obvious scalar inequalities.

The aim of this short survey is to present in an unified  and self-contained way two  recent significant improvement of  the trace inequalities \eqref{VN}-\eqref{Rot} and some of their consequences.

By operator, we mean a linear operator on a finite dimensional Hilbert space. We use interchangeably the terms operator and matrix. Especially, a positive operator means a positive (semi-definite) matrix. Consistently  $\bM_n$ denotes the set of operators on a space of dimension $n$  and  $\bM_n^+$ stands for the  positive part. As many operator inequalities, our results lie in the scope of matrix techniques. Of course, there are versions for operators acting on infinite dimensional, separable Hilbert space (and  operator algebras); we will indicate the slight modifications which might then be necessary.

The rest of this introduction explains why inequalities with unitary orbits are relevant for inequalities involving functional calculus of operators such as the concavity-subadditivity statements  \eqref{VN} and \eqref{Rot}.

That inequalities with unitary orbits  naturally occur can be seen from the following two  elementary facts. Firstly, If $A,B\in \bM_n^+$ are such that $A\ge B$ (that is $A-B$ is positive semi-definite) then, whenever $p>1$,
it does not follow in general that $A^p\ge B^p$. However, for any non-decreasing function $f(t)$, the eigenvalues (arranged in decreasing order and counted with their multiplicities) of $f(A)$ are greater or equal to the corresponding ones of $f(B)$. By the min-max characterization of eigenvalues, this is equivalent to
 \begin{equation}\label{fact}
f(A) \ge Uf(B)U^*
\end{equation}
 for some unitary $U\in\bM_n$. Secondly, if  $A\in\bM_n^+$ and $C\in\bM_n$ is a contraction, then we have $ C^*AC \le UAU^*$ for some unitary  $U\in\bM_n$, i.e., the eigenvalues of $ C^*AC$ are smaller or equal to those of $A$. Note also that $C^*AC=VA^{1/2}CC^*A^{1/2}\le A$ for some unitary $V$, since $TT^*$ and $T^*T$ are unitarily congruent for any operator $T$. The reading of Sections 2 and 3 does not require more knowledge about matrices, see \cite{Bh} for a good background.

The most well-known matrix inequality involving unitary orbits is undoubtedly the triangle inequality due to Thompson \cite{T}: If $X$ and $Y$ are two operators in $\bM_n$, then
\begin{equation}\label{Thompson}
|X+Y| \le U|X|U^* + V|Y|V^*
\end{equation}
for some unitary $U, V\in \bM_n$. Here $|X|:=(X^*X)^{1/2}$ is the positive part of $X$ occurring in the polar decomposition $X=V|X|$ for some unitary $V$. By letting 
$$
X=\begin{pmatrix} A^{1/2}& 0 \\ 0&0 \end{pmatrix}, \qquad Y=\begin{pmatrix} 0& 0\\ B^{1/2}&0 \end{pmatrix}
$$
where $A, B\in\bM_n^+$, the triangle inequality \eqref{Thompson} yields
$
\sqrt{A+B} \le K\sqrt{A}K^* + L\sqrt{B}L^*
$
for some contractions $K, L\in \bM_n$. Thus, for some unitaries 
$U, V\in \bM_n$
\begin{equation}\label{sqrt}
\sqrt{A+B} \le U\sqrt{A}U^* + V\sqrt{B}V^*.
\end{equation}
This inequality for the  function $\sqrt{t}$ is a special case of the main theorem of Section 3.
 
If $f(t)$ is convex on $[0,\infty)$, then \eqref{Rot} is obviously reversed. In case of $f(t)=t^p$ with exponents $p\in[1,2]$ a much stronger inequality holds,
\begin{equation}\label{opconv}
\left(\frac{ A+B}{2}\right)^p \le \frac{A^p +B^p}{2},
\end{equation}
this says that $t^p$ is operator convex for $p\in[1,2]$, and this is not longer true if $p>2$. However by making use of \eqref{fact} and \eqref{opconv} we get, for any $p>1$,
\begin{equation}\label{F-I.1}
\left(\frac{ A+B}{2}\right)^p \le U\frac{A^p +B^p}{2}U^*,
\end{equation}
for some unitary $U\in\bM_n$. This may serve as a motivation for Section 2. It is worthwhile to notice that, in contrast with the theory of operator convex functions, our methods are rather elementary.

\section{ A matrix Jensen type inequality }

\subsection{Jensen type inequalities via unitary orbits}

\noindent In this section we present some extension of (1.1). The most general one involves a unital positive linear map. A linear map $\Phi: \bM_n\to \bM_d$ is unital if $\Phi(I)=I$ where $I$ stands for the identity of any order, and $\Phi$ is positive if $\Phi(A)\in\bM_d^+$ for all $A\in\bM_n^+$. The simplest case is given when $d=1$ by the map
\begin{equation}\label{F-2.1}
A \mapsto \langle h, Ah\rangle
\end{equation}
for some unit vector $h$ (our inner product is linear in the second variable). Restricting this map to the diagonal part (more generally, to any commutative $*$-subalgebra) of $\bM_n$, we have 
\begin{equation}\label{F-2.2}
A \mapsto \langle h, Ah\rangle =\sum_{i=1}^n w_i \lambda_i(A)
\end{equation}
where the $\lambda_i(A)$'s are the eigenvalues of the normal operator $A$ and the $w_i$'s form a probability weight. For this reason,  unital positive linear maps are regarded as non-commutative versions of expectations. If $A$ is Hermitian, and $f(t)$ is a convex function defined on the real line, the Jensen's inequality may be written in term of the map $\Phi$ in \eqref{F-2.1}-\eqref{F-2.2} as 
\begin{equation}\label{Jensen}
f( \langle h, Ah\rangle) \le  \langle h, f(A)h\rangle.
\end{equation}

 The map (2.1) is a special case of a compression. Given an $n$-dimensional Hilbert space ${\mathcal{H}}$ and a $d$-dimensional subspace ${\mathcal{S}} \subset\mathcal{H}$, we have a natural map from the algebra ${\mathrm{L}}(\mathcal{H})$ of operators on $\mathcal{H}$ onto the algebra ${\mathrm{L}}({\mathcal{S}})$, the compression map onto ${\mathcal{S}}$, 
$$
A\mapsto A_{\mathcal{S}} := EA_{|{\mathcal{S}}}, \qquad A\in {\mathrm{L}}(\mathcal{H}),
$$
where $E$ denotes the ortho-projection onto ${\mathcal{S}}$. Identifying ${\mathrm{L}}(\mathcal{H})$  with $\bM_n$ by picking an orthonormal basis of ${\mathcal{H}}$ and ${\mathrm{L}}(\mathcal{S})$ with  $\bM_d$ via an orthonormal basis ${\mathcal{S}}$, we may consider compressions as unital positive linear maps acting from $\bM_n$ onto $\bM_d$, and they are then represented as 
$$
A\mapsto J^*AJ, \qquad A\in\bM_n,
$$
where $J$ is any $n$-by-$d$ matrix such that $J^*J=I$, the identity of order $d$. 

 In view of \eqref{Jensen} it is quite natural to compare for a convex function $f(A_{\mathcal S})$ and $f(A)_{\mathcal S}$ when $A$ is a Hermitian on $\mathcal{H}$, i.e, a Hermitian in $\bM_n$. In this setting, the Jensen inequality \eqref{Jensen} is adapted by using unitary orbits on $\mathcal{S}$. This is actually true for any unital positive linear maps, as stated in Theorem 2.1 below. This is the main result of this section. 
The notation  $\bM_n\{\Omega\}$  stands for the Hermitian part of $\bM_n$ with spectra in an interval $\Omega$ of the real line. 
\vskip 10pt
\begin{theorem}\label{T-2.1} Let $\Phi : \bM_n\to \bM_d$ be a unital positive linear map, let $f(t)$ be a convex function on an interval $\Omega$, and let $A,B\in \bM_n\{\Omega\}$. Then, for some unitary $U,\,V\in\bM_d$,
\begin{equation*}
f(\Phi(A))\le \frac{U\Phi( f(A))U^*+V\Phi( f(A))V^*}{2}.
\end{equation*}
If furthermore $f(t)$ is monotone, then we can take $U=V$. The inequality reverses for concave functions.
\end{theorem}

The next corollaries list some consequences of the theorem.
This  statement for positive linear maps contains several Jensen type inequalities. The simplest one is obtained by taking $\Phi : \bM_{2n}\to \bM_n$,
\begin{equation*}
\Phi\left(\begin{bmatrix}  A &X \\ Y &B\end{bmatrix}\right) :=\frac{A+B}{2}.
\end{equation*}
With $X=Y=0$, Theorem \ref{T-2.1} then says: 

\vskip 5pt
\begin{cor} \label{C-2.2} If  $A,B\in \bM_n\{\Omega\}$ and   $f(t)$ is a convex function on an interval $\Omega$, then, for some unitaries $U,\,V\in\bM_n$,
\begin{equation*}
f\left(\frac{A+B}{2}\right)\le \frac{1}{2}\left\{U\frac{f(A)+f(B)}{2}U^*+V\frac{f(A)+f(B)}{2}V^*\right\}.
\end{equation*}
If furthermore $f(t)$ is monotone, then we can take $U=V$.
\end{cor}

\vskip 5pt
From this corollary we can get a generalization of the famous  Minkowski  inequality,
\begin{equation}\label{Minkowski}
{\det}^{1/n} (A+B) \ge {\det}^{1/n} A + {\det}^{1/n} B, \qquad A,\,B\in \bM_n^+.
\end{equation}
 Equivalently, this says that the Minkowski functional $X\mapsto  \det^{1/n}X$ is concave on the positive cone $\bM_n^+$.
Combined with the concave version of Corollary \ref{C-2.2}, this concavity aspect of \eqref{Minkowski} is  improved as:

\vskip 5pt
\begin{cor}
If $f(t)$ is a non-negative concave function on an interval   $\Omega$ and if $A,B\in \bM_n\{\Omega\}$, then,
\begin{equation*}
{\det}^{1/n} f\left(\frac{A+B}{2}\right) \ge \frac{ {\det}^{1/n} f(A) + {\det}^{1/n} f(B)}{2}.
\end{equation*}
\end{cor}

\vskip 5pt
Corollary \ref{C-2.2}  deals with the simplest convex combination,  the arithmetic mean of two operators. Similar statements holds for weighted  means of several operators. In fact these means may even have operator weights (called $C^*$-convex combinations). An $m$-tuple $\{Z_i\}_{i=1}^m$ in $\bM_n$ is an isometric column if $\sum_{i=1}^m Z_i^*Z_i =I$. We may then perform the $C^*$-convex combination  $\sum_{i=1}^m Z^*_i A_iZ_i $. If all the $A_i$'s
are Hermitian operators in $\bM_n\{\Omega\}$ for some interval $\Omega$, then  so is $\sum_{i=1}^m Z^*_i A_iZ_i $.   Hence, Corollary 2.2 is a very special case of the next one.

\vskip 5pt
\begin{cor}\label{C-2.3} Let $\{Z_i\}_{i=1}^m$ be an isometric column  in  $\bM_n$, let $\{A_i\}_{i=1}^m$ be   in  $\bM_n\{\Omega\}$ and let    $f(t)$ be  a convex function on  $\Omega$. Then, for some unitary $U,\,V\in\bM_n$,
\begin{equation*}
f\left( \sum_{i=1}^m Z^*_i A_iZ_i \right) \le \frac{1}{2} \left\{ U\left(\sum_{i=1}^m Z^*_i f(A_i)Z_i \right)U^*+V\left(\sum_{i=1}^m Z^*_i f(A_i)Z_i\right) V^*\right\}.
\end{equation*}
If furthermore $f(t)$ is monotone, then we can take $U=V$. The inequality reverses for concave functions.
\end{cor}

\vskip 5pt
If all the $A_i$'s are zero except the first one, we obtain an inequality involving a   congruence $Z_1^*A_1Z_1$ with a contraction $Z_1$ (that is $Z_1^*Z_1 \le I$). We state the concave version in the next corollary. It is a matrix version of the basic inequality $f(za) \ge zf(a)$ for a concave function with $f(0)\ge 0$ and real numbers $z,a$ with $z\in[0,1]$.

\vskip 5pt
\begin{cor}\label{C-2.5} Let   $f(t)$ be  a concave function on an interval  $\Omega$ with $0\in\Omega$ and $f(0)\ge 0$, let $A\in \bM_n\{\Omega\}$ and let $Z$ be a contraction in $\bM_n$. Then, for some unitaries $U,\,V\in\bM_n$,
\begin{equation*}
f\left(Z^*AZ\right)\ge \frac{ U\left( Z^*f(A)Z \right)U^*+V\left(Z^* f(A)Z\right) V^*}{2}.
\end{equation*}
If furthermore $f(t)$ is monotone, then we can take $U=V$.
\end{cor}

 For a sub-unital positive linear map $\Phi$, i.e., $\Phi(I)\le I$, it is easy to see that  Theorem 2.1 can be extended in the convex case when $f(0)\le 0$, and in the concave case, when $f(0)\ge 0$ (this sub-unital version is proved in the proof of Corollary 2.7 below). This also contains Corollary \ref{C-2.5}. The above results contains some inequalities for various norms and functionals, as noted in some of the corollaries and remarks below. For instance we have the following Jensen trace inequalities.

\vskip 5pt
\begin{cor}\label{C-2.6} Let    $f(t)$ be  a convex function defined on an interval  $\Omega$, let $\{A_i\}_{i=1}^m$ be   in  $\bM_n\{\Omega\}$, and let $\{Z_i\}_{i=1}^m$ be an isometric column  in  $\bM_n$. Then,
\begin{equation}\label{HP}
{\mathrm{Tr\,}} f\left( \sum_{i=1}^m Z^*_i A_iZ_i \right)\le {\mathrm{Tr\,}} \sum_{i=1}^m Z^*_i f(A)_iZ_i.
\end{equation}
If further $0\in\Omega$ and $f(0)\le 0$, we also have
\begin{equation}\label{BK}
{\mathrm{Tr\,}} f\left( Z^*_1 A_1Z_1 \right)\le {\mathrm{Tr\,}} Z^*_1 f(A_1)Z_1.
\end{equation}
\end{cor}

\vskip 5pt
A typical example of positive linear map on $\bM_n$ is  the Schur multiplication map $A\mapsto Z\circ A$ with an operator $Z\in\bM_n^+$. Here $Z\circ A$ is the entrywise product of $A$ and $Z$. The fact that the Schur multiplication with $Z\in\bM_n^+$  is a positive linear map can be easily checked by restricting the Schur product to positive rank ones operators.
Hence, Theorem 2.1 contains results for the Schur product. In particular, the sub-unital version yields:

\vskip 5pt
\begin{cor}\label{C-2.7} Let   $f(t)$ be  a concave function on   an interval $\Omega$ with $0\in\Omega$ and $f(0)\ge 0$, and let $A\in \bM_n\{\Omega\}$. If $Z\in\bM_n^+$  has  diagonal entries all less than or equal to 1,  then, for some unitaries $U,\,V\in\bM_n$,
\begin{equation*}
f\left(Z\circ A\right)\ge \frac{ U\left( Z\circ f(A) \right)U^*+V\left(Z\circ f(A)\right) V^*}{2}.
\end{equation*}
If furthermore $f(t)$ is monotone, then we can take $U=V$.
\end{cor}

\vskip 5pt
\begin{proof} Let $\Psi:\bM_n\to M_d$ be a positive linear map and suppose that $\Psi$ is sub-unital, i.e., $\Psi(I)=C$ for some contraction $C\in\bM_d^+$. Then the map $\Phi:\bM_{n+1}\to M_d$,
$$
\begin{bmatrix}
A & \vdots \\
\hdots & b
\end{bmatrix}
\mapsto
\Psi(A) + b(I-C)
$$
is  unital. Thus, by Theorem  2.1, If $A\in\bM_n\{{\Omega}\}$ where $\Omega$ contains $0$ and if $f(t)$ is concave on $\Omega$,
$$
f\left(\Phi
(A \oplus 0)
\right) \ge
\frac{U\Phi(f(A\oplus 0))U^*+ V\Phi(f(A\oplus 0))V^*}{2}
$$
for some unitary $U,V\in\bM_d$,
equivalently,
$$
f(\Psi(A)) \ge
\frac{U\{\Psi(A)+f(0)(I-C)\}U^*+ V\{\Psi(A)+f(0)(I-C)\}V^*}{2}
$$
 hence, if further $f(0)\ge 0$, the sub-unital form of Theorem 2.1:
$$
f(\Psi(A)) \ge
\frac{U\Psi(f(A))U^*+ V\Psi(f(A))V^*}{2}.
$$
Applying this to the sub-unital map $\Psi : A\mapsto Z\circ A$ yields the corollary.  \end{proof}

\vskip 5pt Corollary \ref{C-2.7} obviously contains a trace inequality companion to \eqref{BK}. By making use of \eqref{Minkowski} we also have the next determinantal inequality.

\vskip 5pt
\begin{cor}\label{C-2.8} Let   $f(t)$ be  a non-negative concave function on   an interval $\Omega$, $0\in\Omega$, and let $A\in \bM_n\{\Omega\}$. If $Z\in\bM_n^+$  has  diagonal entries all less than or equal to 1,  then, 
\begin{equation}
\det f\left(Z\circ A\right)\ge \det Z\circ f(A).
\end{equation}
\end{cor}

\vskip 5pt
Some other consequences of Theorem 2.1 are given in Subsection 2.2 below, as well as references and related results.

We turn to the proof of Theorem 2.1. Thanks to the next lemma, we will see that it is enough to prove Theorem 2.1 for compressions. By an abelian $*$-subalgebra $\mathcal{A}$ of $\bM_m$ we mean a subalgebra containing the identity of $\bM_m$ and closed under the involution $A\mapsto A^*$. Any  abelian $*$-subalgebra $ {\mathcal{A}}$ of $\bM_m$ is spanned by a total family of ortho-projections, i.e., a family of   mutually orthogonal projections   adding up to the identity. A representation $\pi : \mathcal{A}\to \bM_n$ is a unital linear map such  $\pi(A^*B)=\pi^*(A)\pi(B)$.

\vskip 5pt
 \begin{lemma}  Let $\Phi$ be a unital  positive  map from
 an abelian $*$-subalgebra ${\mathcal{A}}$ of $\bM_n$
to the algebra $\bM_m$ identified as ${\mathrm{L}}({\mathcal{S}})$.  Then, there exists a space ${\mathcal{H}}\supset{\mathcal{S}}$, $\dim {\mathcal{H}}\le nm$, and a
representation
$\pi$ from ${\mathcal{A}}$ to ${\mathrm{L}}({\mathcal{H}})$
such that
$$
\Phi (X)=(\pi(X))_{\mathcal{S}}.
$$
\end{lemma}

 \vskip 5pt
\begin{proof} ${\mathcal{A}}$ is generated by a total family of $k$
projections $E_i$, $i=1,\dots ,k$ (say $E_i$ are rank one, that is
$k=n$). Let $A_i=\Phi(E_i)$, $i=1,\dots ,n$. Since $\sum_{i=1}^n A_i$ is
the identity on ${\mathcal{S}}$, we can find operators $X_{i,j}$ such
that
$$
V=
\begin{pmatrix}
A_1^{1/2}&\dots &A_n^{1/2} \\
X_{1,1} &\dots &X_{n,1} \\
\vdots &\ddots &\vdots \\
X_{1,n-1} &\dots &X_{n,n-1}
\end{pmatrix}
$$
is a unitary operator on ${\mathcal{F}}=\oplus^n{\mathcal{S}}$. Let $R_i$ be the block matrix with the same $i$-th column than $V$ and with all other entries $0$. Then, setting $P_i=R_iR_i^*$, we obtain a total family of projections on ${\mathcal{F}}$ satifying $A_i=(P_i)_{\mathcal{S}}$. We define $\pi$ by $\pi(E_i)=P_i$.  \end{proof}

In the following proof of Theorem 2.1, and in the rest of the paper, 
the eigenvalues of a Hermitian $X$ on an $n$-dimensional space are denoted in non-increasing order as $\lambda_1(X)\ge \cdots\ge\lambda_n(X)$.

\vskip 5pt\noindent 
\begin{proof}  
We consider the convex case.  We first deal with a compression map. Hence $\bM_n$ is identified with ${\mathrm{L}}({\mathcal{H}})$ and 
$
\Phi(A) = A_{\mathcal S}
$
where ${\mathcal{S}}$  is a subspace of ${\mathcal{H}}$.
We may find spectral subspaces ${\mathcal{S}}'$ and ${\mathcal{S}}''$ for $A_{\mathcal{S}}$ and a real $r$ such that

\begin{itemize}
\item[(a)]
 ${\mathcal{S}}={\mathcal{S}}'\oplus{\mathcal{S}}''$,
 
\item[(b)] the spectrum of $A_{\mathcal{S'}}$ lies on $(-\infty,r]$ and the spectrum of  $A_{\mathcal{S''}}$ lies on $[r,\infty)$,
 
 \item[(c)] $f$ is monotone both on $(-\infty,r]\cap\Omega$ and  $[r,\infty)\cap\Omega$.
\end{itemize}

Let $k$ be an integer, $1\le k\le \dim{\mathcal{S}}'$.  There exists a spectral subspace ${\mathcal{F}}\subset{\mathcal{S}}'$ for $A_{\mathcal{S'}}$ (hence for $f(A_{\mathcal{S'}})$), $\dim{\mathcal{F}}=k$, such that
\begin{align*} \lambda_k[f(A_{\mathcal{S'}})] &=\min_{h\in{\mathcal{F}};\ \Vert h\Vert=1} \langle h,f(A_{\mathcal{F}})h \rangle  \\
&= \min\{f(\lambda_1(A_{\mathcal{F}}))\,;\,f(\lambda_k(A_{\mathcal{F}}))\} \\
&= \min_{h\in{\mathcal{F}};\ \Vert h\Vert=1} f(\langle h,A_{\mathcal{F}}h \rangle)  \\
&= \min_{h\in{\mathcal{F}};\ \Vert h\Vert=1} f(\langle h,Ah \rangle) 
\end{align*}
where at the second and third steps we use the monotony of $f$ on $(-\infty,r]$ and the fact that $A_{\mathcal{F}}$'s spectrum lies on $(-\infty,r]$. The convexity of $f$ implies 
$$
f(\langle h,Ah \rangle) \le \langle h,f(A)h \rangle
$$
for all normalized vectors $h$. Therefore, by the minmax principle,
\begin{align*}
\lambda_k[f(A_{\mathcal{S'}})] &\le \min_{h\in{\mathcal{F}};\ \Vert h\Vert=1} \langle h,f(A)h \rangle \\
&\le \lambda_k[f(A)_{\mathcal{S'}}].
\end{align*} 
This  statement is equivalent (by unitary congruence to diagonal matrices) to the existence of a unitary operator $U_0$ on ${\mathcal{S}}'$ such that 
$$
f(A_{\mathcal{S'}})\le U_0 f(A)_{\mathcal{S'}}U_0^*.
$$
(Note that the monotone case is established.) Similarly we get a unitary $V_0$ on ${\mathcal{S}}''$ such that
$$
f(A_{\mathcal{S''}})\le V_0 f(A)_{\mathcal{S''}}V_0^*.
$$
Thus we have
$$
f(A_{\mathcal{S}})\le
\begin{pmatrix} 
U_0 &0 \\ 0&V_0 
\end{pmatrix}  
 \begin{pmatrix} 
f(A)_{\mathcal{S'}} &0 \\ 0& f(A)_{\mathcal{S''}}
\end{pmatrix}  
 \begin{pmatrix} 
U_0^* &0 \\ 0&V_0^* 
\end{pmatrix}. 
$$
Besides we note that, still in respect with the decomposition ${\mathcal{S}}={\mathcal{S}}'\oplus{\mathcal{S}}''$,
$$
\begin{pmatrix} 
f(A)_{\mathcal{S'}} &0 \\ 0& f(A)_{\mathcal{S''}} 
\end{pmatrix}
= \frac{1}{2}\left\{  
\begin{pmatrix} 
I &0 \\ 0& I 
\end{pmatrix}
f(A)_{\mathcal{S}}
\begin{pmatrix} 
I &0 \\ 0& I 
\end{pmatrix}
+
\begin{pmatrix} 
I &0 \\ 0& -I 
\end{pmatrix}
f(A)_{\mathcal{S}}
\begin{pmatrix} 
I &0 \\ 0& -I 
\end{pmatrix}
\right\}.
$$ 
So, letting
$$
U=\begin{pmatrix} 
U_0 &0 \\ 0& V_0 
\end{pmatrix}
\quad{\rm and}\quad
V=
\begin{pmatrix} 
U_0 &0 \\ 0& -V_0 
\end{pmatrix}
$$
we get
\begin{equation}
f(A_{\mathcal{S}}) \le \frac{Uf(A)_{\mathcal{S}}U^* + Vf(A)_{\mathcal{S}}V^*}{2}
\end{equation}
for some unitary $U,V\in{\mathrm{L}}({\mathcal{S}})$, with $U=V$ if $f(t)$ is convex and monotone. This proves the case of compression maps.

Next we turn to the case of a general unital linear map $\Phi:\bM_n\to \bM_m$.
 Let $\mathcal{A}$ be the abelian $*$-subalgebra of $\bM_n$ spanned by $A$. By restricting $\Phi$ to $\mathcal{A}$ and by identifying $\bM_m$ with ${\mathrm{L}}({\mathcal{S}})$, Lemma 2.9 shows that $\Phi(X)=(\pi(X))_{\mathcal{S}}$ for all $X\in\mathcal{A}$.
Since $f$ and $\pi$ commutes, $f(\pi(A))=\pi(f(A))$, we have from the compression case some unitary $U, V\in {\mathrm{L}}({\mathcal{S}})=\bM_m$ such that,
\begin{align*}
f(\Phi(A) ) &=f((\pi(A))_{\mathcal{S}})  \\ &\le \frac{U(f(\pi(A))_{\mathcal{S}}U^* + V(f(\pi(A))_{\mathcal{S}}V^*}{2} \\
&= \frac{U(\pi(f(A))_{\mathcal{S}}U^* + V(\pi (f(A))_{\mathcal{S}}V^*}{2} \\
&= \frac{U\Phi(f(A))U^* + V\Phi(f(A))V^*}{2},
\end{align*}
where we can take $U=V$ if the function is convex and monotone.
\end{proof}

\vskip 5pt
The following is an application of Theorem 2.1 to norm inequalities.
A norm $\|\cdot\|$ on $\bM_n$ is a symmetric norm if $\|A\|=\|UAV\|$ for all $A\in\bM_n$ and all unitary $U,V\in\bM_n$. These norms are also called unitarily invariant norms. They contain the Schatten $p$-norms. The polar decomposition shows that a symmetric norm is well defined by its value on the positive cone $\bM_n^+$. The map on $\bM_n^+$, $A\mapsto\|A\|$ is invariant under unitary congruence and is  subadditive. There are also some interesting, related superadditive functionals.
Fix $p<0$. The map $X\mapsto\|X^p\|$ is  continuous on the invertible part of $\bM_n^+$. If $X\in\bM_n^+$ is not invertible,  setting $\| X^p\|:=0$, we obtain a continuous map on $\bM_n^+$. 

\vskip 5pt
\begin{cor} Let $A,B\in\bM_n^+$ and let $p<0$. Then, for all symmetric norms,
\begin{equation} \label{derived}
\|\, (A+B)^p\|^{1/p} \,\ge\, \| A^p\|^{1/p}+ \| B^p\|^{1/p}.
\end{equation}
\end{cor}

\vskip 5pt
\begin{proof}
 We  will apply Theorem 2.1  to the monotone convex function on  $(0,\infty)$, $t\mapsto t^p$. 
First, assume that $A,B\in\bM_n^+$ are such that $\|A^p\|=\|B^p\|=1$ and let $s\in[0,1]$. Then, thanks to Theorem 2.1 (or Corollary 2.4),
$$
\|(sA+(1-s)B)^p\| \le \| sA^p + (1-s)B^p \| 
\le  s\|A^p\| + (1-s)\|B^p\| 
=1,
$$
hence
\begin{equation}\label{F-2.14}
\|(sA+(1-s)B)^p\|^{1/p} \ge 1.
\end{equation}
Now, for general invertible $A,B\in\bM_n^+$, insert $A/\|A^p\|^{1/p}$ and $B/\|B^p\|^{1/p}$ in place of $A, B$ in \eqref{F-2.14} and take
$$
s=\frac{\|A^p\|^{1/p}}{\|A^p\|^{1/p} +\|B^p\|^{1/p}}.
$$
This yields \eqref{derived}.
\end{proof}

\vskip 5pt
Corollary 2.10 implies Minkowsky's determinantal inequality (2.4). Indeed, in (2.9) take the  norm on $\bM_n^+$ defined by $\| A\|:=\frac{1}{n}{\mathrm{Tr\,}}A$,
and note that $\det^{1/n} A=\lim_{p\nearrow 0} \| A^p\|^{1/p}$.
Hence,  the superadditivity of $A\mapsto \|A^p\|^{1/p}$ for $p<0$ entails the superadditivity of $A\mapsto \det^{1/n} A$.

If we apply Theorem 2.1 (or Corollary 2.2) to the convex function on the real line $t\mapsto |t|$ we obtain: {\it If $A,B\in\bM_n$ are Hermitian, then
\begin{equation*}
|A+B| \le \frac{U(|A|+|B|)U^*+V(|A|+|B|)V^*}{2}
\end{equation*}
for some unitaries $U,V\in\bM_n$.} In fact, we can take $U=I$ and this remains true for normal operators $A,B$. This is shown in the proof of the following proposition.

\vskip 5pt
\begin{prop} If $f(t)$ is a nondecreasing convex function on $[0,\infty)$ and if $Z\in\bM_n$ has a Cartesian decomposition $Z=A+iB$, then, for some unitaries $U,V\in\bM_n$,
\begin{equation*}
f(|Z|) \le \frac{Uf(|A|+|B|)U^*+Vf(|A|+|B|)V^*}{2}.
\end{equation*}
\end{prop}

\vskip 5pt
\begin{proof} let $X$, $Y$ be two normal operators in $\bM_n$. Then, the following operators in $\bM_{2n}$ are positive semi-definite,
$$
\begin{pmatrix}
|X| &X^* \\
X&|X|
\end{pmatrix} \ge 0,
\qquad 
\begin{pmatrix}
|Y| &Y^* \\
Y&|Y|
\end{pmatrix} \ge 0,
$$
and consequently
$$
\begin{pmatrix}
|X| +|Y|&X^*+Y^* \\
X+Y&|X|+|Y| 
\end{pmatrix}
\ge 0.
$$
Next, let $W$ be the unitary part in the polar decomposition $X+Y=W|X+Y|$. Then
$$
\begin{pmatrix}
I&-W^* 
\end{pmatrix}
\begin{pmatrix}
|X| +|Y|&X^*+Y^* \\
X+Y&|X|+|Y| 
\end{pmatrix}
\begin{pmatrix}
I \\ -W 
\end{pmatrix}
\ge0,
$$
that is 
$$
|X|+|Y|+W^*(|X|+|Y|)W- 2|X+Y|\ge 0.
$$
Equivalently,
\begin{equation}\label{F-2.11}
|X+Y| \le \frac{|X|+|Y|+W^*(|X|+|Y|)W}{2}.
\end{equation}
Letting $X=A$ and $Y=iB$, and applying $f(t)$ to both sides of \eqref{F-2.11}, 
  Corollary 2.2 completes the proof since $f(t)$ is nondecreasing and convex.
\end{proof}

\vskip 5pt
\begin{prop} If $f(t)$ is a nondecreasing convex function on $[0,\infty)$ and if $A,B\in\bM_n$ are Hermitian, then, for some unitaries $U,V\in\bM_n$,
\begin{equation*}
f((A+B)_+) \le \frac{Uf(A_++B_+)U^*+Vf(A_++B_+)V^*}{2}.
\end{equation*}
\end{prop}

\vskip 5pt
\begin{proof} Here $A_+:=(A+|A|)/2$.
Note that $A+B\le  A_++ B_+$. Let $E$ be the projection onto ${\mathrm{ran\,}}(A+B)_+$ and let $F$ be the projection onto $\ker(A+B)_+$
Since
$(A+B)_+ =E(A+B)E$, we have
$$
(A+B)_+\le  E(A_++ B_+)E + F(A_++B_+)F,
$$
equivalently
\begin{equation}
(A+B)_+\le \frac{(A_++ B_+) +W(A_++ B_+)W^*}{2}
\end{equation}
where $W= E-F$ is a unitary. Applying Corollary 2.2 completes the proof.
 \end{proof}

\subsection{Comments and references}

\noindent
In this second part of Section  2, we collect few remarks which complete Theorem 2.1 and the above corollaries. Good references for  positive maps and operator convex functions  are the nice  survey   and book \cite{Hiai1} and \cite{Bh}.

\vskip 10pt
\begin{remark} Theorem 2.1  appears in \cite{B2}. It is stated therein for compressions maps and for the case of  or $*$-convex combinations given in Corollaries 2.4 and 2.5 (the monotone case was earlier  obtained in \cite{B1}).  That the compression case immediately entails the general case of an arbitrary unital positive map is mentioned in some subsequent papers, for instance in \cite{AB} where some  inequalities for Schur products are pointed out. From the Choi-Kraus representation of completely positive linear maps,  readers with a background on positive maps may also notice that Corollary 2.5 and Theorem 2.1 are equivalent.
For scalar convex combinations and with the assumption that $f(t)$ is non-decreasing, Theorem 2.1  is first noted in Brown-Kosaki's paper \cite{BK}; with these assumptions, it is also obtained in Aujla-Silva's paper \cite{AS}.
\end{remark}

\vskip 10pt
\begin{remark} Let $g(t)$ denote either the convex function $t\mapsto |t|$ or $t\mapsto t_+$.  Let $A,B\in\bM_n$ be Hermitian. Then  (2.11) and (2.12) show that
$$
g\left(\frac{A+B}{2}\right) -\frac{g(A)+g(B)}{4} \le V\frac{g(A)+g(B)}{4}V^*
$$
for some unitary $V\in\bM_n$. It would be interesting to characterize convex functions for which such a relation holds.
\end{remark}

\vskip 10pt
\begin{remark} Theorem 2.1 holds for operators acting on infinite dimensional spaces, with an additional $rI$ term. We state here the monotone version. ${\mathcal{H}}$ and ${\mathcal{S}}$ are two separable Hilbert spaces  and $r>0$ is fixed. {\it
 Let $\Phi : {\mathrm{L}}({\mathcal{H}})\to  {\mathrm{L}}({\mathcal{H}})$ be a unital positive linear map, let  $f(t)$ be a monotone convex function on $(-\infty,\infty)$ and let
 $A,B\in  {\mathrm{L}}({\mathcal{H}}) $  be Hermitian. Then, for some unitary $U\in {\mathrm{L}}({\mathcal{S}})$,}
\begin{equation}\label{F-2.13}
f(\Phi(A))\le U\Phi( f(A))U^* +rI.
\end{equation}
The proof is given in the first author's thesis when $\Phi$ is a compression map, this entail the general case. For convenience, the proof is given at the end of this section.
\end{remark}

\vskip 5pt
\begin{remark}
The trace inequality \eqref{HP} is due to Hansen-Pedersen \cite{HP2}, and the special case \eqref{BK} is due to Brown-Kosaki. Corollaries 2.3 considerably improves \eqref{HP}: In case of a monotony assumption on the convex function $f(t)$, we have eigenvalue inequalities; and, in the general case we may still infer the majorization relation
\begin{equation}\label{F-2.9}
\sigma_k \left[f\left( \sum_{i=1}^m Z^*_i A_iZ_i \right)\right]\le \sigma_k\left[ \sum_{i=1}^m Z^*_i f(A)_iZ_i \right], \qquad k=1,\ldots,n,
\end{equation}
where $\sigma_k[X]:=\sum_{j=1}^k\lambda_j[X]$ is the sum of the $k$  largest eigenvalues of a Hermitian $X$. In fact,  the basic relation $\sigma_k[X]=\max {\mathrm{Tr\,}} XE$, where the maximum runs over all rank $k$ projections $E$, shows that $\sigma_k[\cdot]$ is convex, increasing on the Hermitian part of $\bM_n$ so that \eqref{F-2.9} is an immediate consequence of Theorem 2.1. The theorem also entails (see \cite{B2} for details) a rather unexpected eigenvalue inequality:
$$
\lambda_{2k-1}\left[ f\left( \sum_{i=1}^m Z^*_i A_iZ_i \right)\right]  \le \lambda_k\left| \sum_{i=1}^m Z^*_i f(A)_iZ_i \right] , \qquad  1\le k\le  (n+1)/2.
$$
\end{remark}

\vskip 5pt
\begin{remark}  Choi's inequality \cite{Choi} claims: {\it for an operator convex function $f(t)$ on  $\Omega$,
\begin{equation}\label{Choi}
f(\Phi(A)) \le \Phi(f(A))
\end{equation}
for all $A\in\bM_n\{\Omega\}$ and all unital positive linear map.} Thus
 Theorem 2.1 is a substitute of Choi's inequality for a general convex function. In the special of a compression map, then \eqref{Choi} is Davis' inequality \cite{Davis}, a famous characterization of operator convexity. The most well-known case of Davis' inequality is for the inverse map on positive definite matrices, it is then an old classical fact of Linear Algebra. Exactly as Theorem 2.1 entails Corollary 2.3, Choi's inequality contains
Hansen-Pedersen's inequality \cite{HP1}, \cite{HP2}: {\it If $f(t)$ is operator convex on $\Omega$, then
\begin{equation*}
f\left( \sum_{i=1}^m Z^*_i A_iZ_i \right) \le 
\sum_{i=1}^m Z^*_i f(A)_iZ_i 
\end{equation*}
for  all  unitary columns $\{Z_i\}_{i=1}^m$ in $\bM_n$ and $A_i\in\bM_n\{\Omega\}$, $i=1,\ldots,m$.}
For  operator concave functions,  the inequality reverses. A special case is Hansen's  inequality  \cite{H}: {\it  if $f(t)$ is operator concave on $\Omega$, $0\in\Omega$ and $f(0)\ge 0$, then
\begin{equation}\label{Hansen}
f(Z^*AZ) \ge Z^*f(A)Z
\end{equation}
for all $A\in\bM_n\{\Omega\}$ and all contractions $Z\in\bM_n$.}
\end{remark}

\vskip 10pt
\begin{remark} Hansen's inequality \eqref{Hansen} may be formulated with an expansive operator $Z\in\bM_n$, i.e., $Z^*Z\ge I$; then \eqref{Hansen} obviously reverses. We might expect that in a similar way, Corollary 2.5 or the Brown-Kosaki trace inequality reverses. But this does not hold.  Corollary 2.5 can not reverse when $Z$ is expansive, even under the monotony assumption on $f(t)$. An unexpected positivity assumption is necessary, and we must confine to weaker inequalities, such as trace inequalities: {\it if $f(t)$ is a concave function on the positive half-line with $f(0) \ge0$, then,
\begin{equation*}
{\mathrm{Tr\,}} f\left( Z^* AZ \right)\le {\mathrm{Tr\,}} Z f(A)Z
\end{equation*}
for all $A\in\bM_n^+$ and all expansive $Z\in\bM_n$.} For a proof, see \cite{B1}   and also \cite{B3}, \cite{BL} where remarkable extensions to norm inequalities are given.
\end{remark}

\vskip 10pt
\begin{remark} Lemma 2.9 is a part of Stinespring' s theory of positive and completely positive linear maps in the influential 1955 paper \cite{S}. The proof given here is somewhat original
and is taken from \cite{AB}. Note that in the curse of the proof, we prove Naimark's dilation theorem: {\it If $\{A_i\}_{i=1}^n$ are positive operators on a space ${\mathcal{S}}$ such that $\sum_{i=1}^n A_i \le I$, then there exist some mutually orthogonal projections $\{P_i\}_{i=1}^n$  on a larger space ${\mathcal{H}}\supset{\mathcal{S}}$ such that $(P_i)_{\mathcal{S}}=A_i$, ($1\le i\le n$). }
\end{remark}

\vskip 10pt
\begin{remark}  Given a symmetric norm $\|\cdot\|$ and $p<0$, the  functionals defined on $\bM_n^+$,
$A\mapsto \|A^p\|^{1/p}$, are introduced in \cite{BH2} and called {\it derived anti-norms}. Corollary 2.10 is given therein,  \cite[Proposition 4.6]{BH2}. The above proof is much simpler than the original one. For more details  and many results on {\it anti-norms} and {\it derived anti-norms}, often in connection with Theorem 2.1, see  \cite{BH1} and \cite{BH2}. Several results in these papers are generalizations of Corollary 2.3. By using (1.7) and arguing as in the proof of Corollary 2.10, we may derive the triangle inequality for the Schatten $p$-norms, i.e.,
$$
\{ {\mathrm{Tr}\,} (A+B)^p \}^{1/p} \le \{ {\mathrm{Tr\,}} A^p \}^{1/p} + \{ {\mathrm{Tr\,}} B^p \}^{1/p}, \qquad A,B\in\bM_n^+, \  p>1.
$$
\end{remark}

\vskip 10pt
\begin{remark}  The inequality  \eqref{F-2.11} for normal operators can be extended to general $A,B\in\bM_n$, with a similar proof, as 
$$
|A+B|\le\frac{|A|+|B|+V(|A^*|+|B^*|)V^*}{2}
$$
for some unitary $V\in\bM_n$. This is pointed out in \cite{BR}. This is still true for operators $A,B$ in a von Neumann algebra $\mathcal{M}$ with $V$ a partial isometry in $\mathcal{M}$. If $\mathcal{M}$ is endowed with a regular trace, this gives a short, simple proof of the triangle inequality for the trace norm on $\mathcal{M}$. 
Inequality \eqref{F-2.11} raises the question of comparison $|A+B|$ and $|A|+|B|$. The following result is given in \cite{Lee2}.
 {\it Let $A_1,\cdots, A_m$ be invertible operators  with
condition numbers dominated by  $\omega > 0 $. Then}
$$
|A_1+\cdots +A_m|\leq \frac{\omega+1}{2\sqrt{\omega}}(|A_1|+\cdots
+|A_m|).
$$
Here the condition number of an invertible operator $A$ on a Hilbert space is $\|A\|\|A^{-1}\|^{-1}$. Note that the bound is independent of the number of operators. Though it is a rather low bound, it is not known whether it is sharp.
\end{remark}

\vskip 5pt
It remains to give a proof of  the infinite dimensional version \eqref{F-2.13} of the monotone case of Theorem 2.1, the non-monotone case following in a similar way to the finite dimensional version. As for finite dimensional spaces,  we may assume that $\Phi$ is a compression map, thus we consider a subspace $\mathcal{S}\subset{\mathcal H}$ and the map $A\mapsto A_{\mathcal{S}}$. By replacing $f(t)$ by $f(-t)$ and $A$ by $-A$, we may also assume that $f(t)$ is nondecreasing.

If $X$ is a Hermitian on $\mathcal{H}$, we define a sequence of numbers $\{\lambda_k(X)\}_{k=1}^{\infty}$,
$$
\lambda_k(X)=\sup_{\{{\cal F}\,:\,\dim{\cal F}=k\}}\,\inf_{\{h\in{\cal F}\,:\,\Vert h\Vert=1\}}\langle h, Xh\rangle
$$ 
where the supremum runs over $k$-dimensional subspaces. Note that $\{\lambda_k(X)\}_{k=1}^{\infty}$  is a non-increasing sequence whose limit is the upper bound of the essential spectrum of $X$. We also define $\{\lambda_{-k}(X)\}_{k=1}^{\infty}$,
$$
\lambda_{-k}(X)=\sup_{\{{\cal F}\,:\,{\rm codim\,}{\cal F}=k-1\}}\,\inf_{\{h\in{\cal F}\,:\,\Vert h\Vert=1\}}\langle h, Xh\rangle.
$$ 
Then, $\{\lambda_{-k}(X)\}_{k=1}^{\infty}$ is a nondecreasing sequence whose limit is the lower bound of the essential spectrum of $X$. The following  fact (a) is obvious and fact (b) is easily checked.
\begin{itemize}
\item[(a)] If  $X\le Y$, then $\lambda_k(X)\le\lambda_k(Y)$ and $\lambda_{-k}(X)\le\lambda_{-k}(Y)$ for all $k=1,\cdots$.

\item[(b)] if $r>0$ and $X,\,Y$ are Hermitian, $\lambda_k(X)\le\lambda_k(Y)$ and $\lambda_{-k}(X)\le\lambda_{-k}(Y)$,for all $k=1,\dots$, then  $X\le UYU^* + rI$ for some unitary $U$.
\end{itemize}
 These facts show that, given $r>0$, two Hermitians $X$, $Y$ with $X\le Y$,  and a continuous nondecreasing  function $\phi$,  there exists a unitary  $U$ such that $\phi(X)\le U\phi(Y)U^*+rI$. 

By fact (b) it suffices to show that 
\begin{equation}\label{F-2.a}
\lambda_k(f(A_{\cal S}))\le\lambda_k(f(A)_{\cal S})
\end{equation}
and 
\begin{equation}\label{F-2.b}
\lambda_{-k}(f(A_{\cal S}))\le\lambda_{-k}(f(A)_{\cal S})
\end{equation} 
for all $k=1,\cdots$. Now, we prove \eqref{F-2.b}  and  distinguish two cases:
\vskip 5pt\noindent
 1.\ $\lambda_{-k}(A_{\cal S})$ is an eigenvalue of $A_{\cal S}$. Then, for $1\le j\le k$, $\lambda_{-j}(f(A_{\cal S}))$ are eigenvalues for $f(A_{\cal S})$. Consequently, there exists a subspace
${\cal F}\subset{\cal S}$, ${\rm codim}_{\cal S}\,{\cal F}=k-1$, such that
\begin{align*}
\lambda_{-k}(f(A_{\cal S})) =& \min_{\{h\in{\cal F}\,:\,\Vert h\Vert=1\}}\langle h, f(A_{\cal S})h\rangle \\
=& \min_{\{h\in{\cal F}\,:\,\Vert h\Vert=1\}}f(\langle h, A_{\cal S}h\rangle) \\
\le& \inf_{\{h\in{\cal F}\,:\,\Vert h\Vert=1\}}\langle h, f(A)h\rangle \le \lambda_{-k}(f(A)_{\cal S})
\end{align*}
where we have used that $f$ is non-decreasing and convex.

\noindent
 2.\ $\lambda_{-k}(A_{\cal S})$ is not an eigenvalue of $A_{\cal S}$ (so, $\lambda_{-k}(A_{\cal S})$ is the lower bound of the essential spectrum of $A_{\cal S}$). 
 Fix $\varepsilon>0$ and choose  $\delta>0$ 
 such that $|f(x)-f(y)|\le \varepsilon$ for all $x$, $y$ are in the convex hull of the spectrum of $A$ with $|x-y|\le\delta$. There exists a subspace ${\cal F}\subset{\cal S}$, ${\rm codim}_{\cal S}\,{\cal F}=k-1$, such that
$$
\lambda_{-k}(A_{\cal S}) \le \inf_{\{h\in{\cal F}\,:\,\Vert h\Vert=1\}}\langle h, A_{\cal S}h\rangle +\delta
.
$$
Since $f$ is continuous nondecreasing we have $f(\lambda_{-k}(A_{\cal S}))=\lambda_{-k}(f(A_{\cal S}))$ so that, as $f$ is nondecreasing,
$$
\lambda_{-k}(f(A_{\cal S})) \le f\left(\inf_{\{h\in{\cal F}\,:\,\Vert h\Vert=1\}}\langle h, A_{\cal S}h\rangle +\delta\right).
$$ 
Consequently,
$$
\lambda_{-k}(f(A_{\cal S})) \le \inf_{\{h\in{\cal F}\,:\,\Vert h\Vert=1\}}f(\langle h, A_{\cal S}h\rangle) + \varepsilon,
$$ 
 so, using the convexity of $f$ and the definition of $\lambda_{-k}(\cdot)$, we get
$$
\lambda_{-k}(f(A_{\cal S}))\le \lambda_{-k}(f(A)_{\cal S}) + \varepsilon.
$$
By letting $\varepsilon\longrightarrow 0$, the proof \eqref{F-2.b} is complete.  The proof of \eqref{F-2.a} is similar.
Thus \eqref{F-2.13} is established.

\vskip 5pt

\section{A matrix subadditivity inequality}

\subsection{Sub/super-additivity inequalities via unitary orbits}

This section deals with some recent subadditive properties for concave functions, and similarly superadditive properties of convex functions. The main result is:

\vskip 10pt
\begin{theorem}\label{T-3.1} Let $f(t)$ be a monotone concave function on $[0,\infty)$ with $f(0)\ge 0$ and let
$A,B\in\bM_n^+$. Then, for some unitaries $U,V\in\bM_n$,
$$
f(A+B)\le Uf(A)U^* + Vf(B)V^*.
$$
\end{theorem}

\vskip 5pt\noindent
Thus ,  the obvious scalar  inequality $f(a+b)\le f(a)+f(b)$   can be extended to positive  matrices $A$ and $B$ by considering element in the unitary orbits of $f(A)$ and $f(B)$.  This inequality via unitary orbits considerably improves the famous Rotfel'd trace inequality \eqref{Rot} for a non-negative concave function on the positive half-line, 
and its symmetric norm version
\begin{equation}
\| f(A+B)\|  \le \|f(A)\| + \|f(B)\|
\end{equation}
for all $A,B\in \bM_n^+$ and all symmetric norms $\|\cdot\|$ on $\bM_n$. 

Of course Theorem \ref{T-3.1} is equivalent to the next statement for  convex functions:

\vskip 10pt
\begin{cor} Let $g(t)$ be a monotone convex function on $[0,\infty)$ with $g(0)\le0$ and let
$A,B\in \bM_n^+$ . Then, for some unitaries $U,V\in\bM_n$,
\begin{equation}\label{supad}
g(A+B)\ge Ug(A)U^* + Vg(B)V^*.
\end{equation}
\end{cor}

 \begin{proof}  It suffices to prove the convex version, Corollary 3.2. We may confine the proof to the case $g(0)=0$ as if \eqref{supad} holds for a function $g(t)$ then it also holds for $g(t)-\alpha$ for any $\alpha>0$. This assumption combined with the monotony of $g(t)$ entails that $g(t)$ has a constant sign $\varepsilon\in\{-1,1\}$, hence $g(t)=\varepsilon |g|(t)$.

We may also assume that $A+B$ is invertible. Then
 $$
 A=X(A+B)X^*
 \quad {\rm and} \quad
  B=Y(A+B)Y^*
 $$
where
 $X=A^{1/2}(A+B)^{-1/2}$ and $Y=B^{1/2}(A+B)^{-1/2}$ are
 contractions. For any $T\in \bM_n,$ $T^*T$ and $TT^*$ are unitarily
congruent. Hence, using Corollary \ref{C-2.5} we have  two unitary operators $U_0$ and $U$ such that
 \begin{align*}
 g(A) &= g(X(A+B)X^*)\\
&\leq U_0 Xg(A+B)X^* U_0^*\\
&= \varepsilon U^*(|g|(A+B))^{1/2}X^*X(|g|(A+B))^{1/2} U,
\end{align*}
so,
\begin{equation}
Ug(A)U^*\leq \varepsilon(|g|(A+B))^{1/2}X^*X(|g|(A+B))^{1/2}.
\end{equation}
Similarly there exists a unitary operator $V$ such that
\begin{equation}
Vg(B)V^*\leq\varepsilon (|g|(A+B))^{1/2}Y^*Y(|g|(A+B))^{1/2}.
\end{equation}
Adding (3.3) and (3.4) we get
$$
 Ug(A)U^*+Vg(B)V^*\leq g (A+B)
$$
since $X^*X+Y^*Y=I_n.$ \end{proof}

\vskip 10pt
The following corollary is matrix version of another obvious scalar inequality.
\vskip 10pt
 \begin{cor}
 Let $f:[0,\infty)\to [0,\infty)$ be concave and let
$A,B\in\bM_n$ be Hermitian.
  Then, for some unitaries $U,V\in\bM_n$,
$$
 Uf(A)U^*-Vf(B)V^*\leq f (|A-B|).
$$
\end{cor}

\vskip 10pt
 \begin{proof} Note that
 $$
 A\leq |A-B|+B.
 $$
 Since $f(t)$ is non-decreasing  and concave there exists  unitaries $W,\,S,\,T$ such that
 $$
 Wf(A)W^*\leq f(|A-B|+B)\leq Sf(|A-B|)S^*+Tf(B)T^*.
 $$
Hence, we have
$$
 Uf(A)U^*-Vf(B)V^*\leq f (|A-B|)
$$
for some unitaries $U,\, V.$ \end{proof}

We can employ Theorem 3.1 to get an elegant inequality for positive block-matrices,
 $$\begin{bmatrix} A &X \\
X^* &B\end{bmatrix}\in \bM_{n+m}^+, \qquad A\in\bM_n^+, \,  B\in\bM_m^+,
$$
which nicely extend (3.1). To this end we need an interesting
  decomposition lemma  for elements in $\bM_{n+m}^+$.
 
\begin{lemma} For every matrix in  $\bM_{n+m}^+$ written in blocks, we have a decomposition
\begin{equation}
\begin{bmatrix} A &X \\
X^* &B\end{bmatrix} = U
\begin{bmatrix} A &0 \\
0 &0\end{bmatrix} U^* +
V\begin{bmatrix} 0 &0 \\
0 &B\end{bmatrix} V^*
\end{equation}
for some unitaries $U,\,V\in  \bM_{n+m}$.
\end{lemma}

\vskip 5pt
\begin{proof} To obtain this decomposition of the positive semi-definite block matrix, factorize it as a square of positive matrices,
\begin{equation*}
\begin{bmatrix} A &X \\
X^* &B\end{bmatrix} =
\begin{bmatrix} C &Y \\
Y^* &D\end{bmatrix}
\begin{bmatrix} C &Y \\
Y^* &D\end{bmatrix}
\end{equation*}
and observe that it can be written as
\begin{equation*}
\begin{bmatrix} C &0 \\
Y^* &0\end{bmatrix}
\begin{bmatrix} C &Y \\
0 &0\end{bmatrix} +
\begin{bmatrix} 0 &Y \\
 0&D\end{bmatrix}
\begin{bmatrix} 0 &0 \\
Y^* &D\end{bmatrix} = T^*T + S^*S.
\end{equation*}
Then, use the fact that $T^*T$ and $S^*S$ are unitarily congruent to
$$
TT^*= \begin{bmatrix} A &0 \\
0 &0\end{bmatrix}
\quad \mathrm{and} 
\quad
SS^*=\begin{bmatrix} 0 &0 \\
0 &B\end{bmatrix},
$$
completing the proof of the decomposition.
\end{proof}

Combined with Theorem \ref{T-3.1}, the lemma yields a norm inequality for block-matrices. A symmetric norm on $\bM_{n+m}$ induces a symmetric norm on  $\bM_{n}$, via $\|A\|=\|A \oplus 0\|$.

\vskip 10pt
\begin{cor} Let $f(t)$ be a non-negative concave
function on $[0,\infty)$. Then, given an arbitrary partitioned
positive semi-definite matrix,
$$
\left\| \,f\left( \begin{bmatrix} A &X \\ X^* &B\end{bmatrix}\right)
 \right\|
\le \left\| f(A) \right\| +  \left\| f(B) \right\|
$$
for all symmetric norms.
\end{cor}

\vskip 10pt
\begin{proof}
From (3.5) and Theorem 3.1, we have 
\begin{equation*}
f\left(\begin{bmatrix} A &X \\
X^* &B\end{bmatrix}\right) = U
\begin{bmatrix} f(A) &0 \\
0 &f(0)I\end{bmatrix} U^* +
V\begin{bmatrix} f(0)I &0 \\
0 &f(B)\end{bmatrix} V^*
\end{equation*}
for some unitaries $U,\,V\in  \bM_{n+m}$. The result then follows from the simple fact
that symmetric norms are nondecreasing functions of the singular values.
\end{proof}

\vskip 10pt
 Applied to  $X=A^{1/2}B^{1/2}$, this result yields the  Rotfel'd type inequalities (1.2)-(3.1), indeed, 
$$
\begin{bmatrix} A &X \\
X^* &B\end{bmatrix}=\begin{bmatrix} A^{1/2} &0 \\
B^{1/2} &0\end{bmatrix}\begin{bmatrix} A^{1/2} &B^{1/2} \\
0 &0\end{bmatrix}
$$
is then  unitarily equivalent to $(A+B)\oplus 0$.
 In case of the trace norm, the above result may be restated as a trace inequality without any non-negative assumption: {\it For all concave functions $f(t)$ on the positive half-line and all positive block-matrices,}
$$
{\mathrm{Tr\,}} f\left( \begin{bmatrix} A &X \\ X^* &B\end{bmatrix}\right)
\le {\mathrm{Tr\,}} f(A) +  {\mathrm{Tr\,}} f(B).
$$
The case of $f(t)=\log t$ then gives Fisher's inequality, 
$$
\det \begin{bmatrix} A &X \\ X^* &B\end{bmatrix} \le \det A \det B.
$$

Theorem 3.1 may  be used to extend another classical (superadditive and concavity) property of the determinant, the Minkowsky inequality (2.4). We have the following extension:

\begin{cor}  If $g:[0,\infty)\to [0,\infty)$ is a convex function, $g(0)=0$,
and  $A,\,B\in \bM_n^+$, then,
\begin{equation*}
{\det}^{1/n} g(A+B) \ge {\det}^{1/n} g(A) + {\det}^{1/n} g(B).
\end{equation*}
\end{cor}

\vskip 10pt
As another example of combination of Theorem 3.1 and (3.5), we have:

\vskip 10pt
\begin{cor} Let $f:[0,\infty)\to [0,\infty)$ be concave and let $A=(a_{i,j})$
be a positive semi-definite matrix in $\bM_n$. Then, for some
rank one ortho-projections $\{E_i\}_{i=1}^n$ in $\bM_n$,
$$
f(A)\le \sum_{i=1}^n f(a_{i,i})E_i.
$$
\end{cor}

\begin{proof} By a limit argument, we may assume that $A$ is invertible, and hence we may also assume that
$f(0)=0$,
indeed if the spectrum of $A$ lies in an interval $[r,s]$, $r>0$, we may replace $f(t)$
 by  any concave function on
$[0,\infty)$ such that $\tilde{f}(0)=0$ and $\tilde{f}(t)=f(t)$ for
$t\in[r,s]$.   By a repetition of (3.5) we have
$$
A=\sum_{i=1}^n a_{i,i} F_i
$$
for some
rank one ortho-projections $\{F_i\}_{i=1}^n$ in $\bM_n$. An application of Theorem 3.1 yields
$$
f(A) \le \sum_{i=1}^n U_i f(a_{i,i} F_i)U_i^*
$$
for some unitary operators $\{U_i\}_{i=1}^n$. Since $f(0)=0$, for each $i$, $U_i f(a_{i,i} F_i)U_i^*=f(a_{i,i})E_i$ for some rank one projection $E_i$.
\end{proof}

\vskip 5pt
 Corollary 3.7 refines the standard majorization inequality relating a positive semi-definite  $n$-by-$n$ matrix and its diagonal part, $${\mathrm{Tr\,}} f(A) \le  \sum_{i=1}^d f(a_{i,i}). $$

\subsection{Comments and references}

\vskip 5pt
\begin{remark} Theorem 3.1, Corollaries 3.2 and 3.3 are from \cite{AB}. In case of positive operators acting  on an infinite dimensional, separable Hilbert space, we have a version of Theorem 3.1 with an additional $r$I term in the RHS, as in (2.13).
\end{remark}

\vskip 5pt
\begin{remark} The decomposition of a positive block-matrix in Lemma 3.4 is due to the authors. It is used in \cite{Lee1} to obtain the norm inequality stated in Corollary 3.5. The next  two  Corollaries 3.6 and 3.7 are new, though already announced in \cite{BH1}.
\end{remark}

\vskip 5pt
\begin{remark}
The concavity requirement on $f(t)$ in Rotfel'd inequality (1.2) and hence in Theorem 3.1 cannot be relaxed to a mere
superadditivity assumption; indeed take for $s,t>0$,
$$
A=\frac{1}{2}
\begin{bmatrix} s&\sqrt{st} \\
\sqrt{st}&t
\end{bmatrix},
\qquad B=\frac{1}{2}
\begin{bmatrix} s&-\sqrt{st} \\
-\sqrt{st}&t
\end{bmatrix},
$$
and observe that the trace inequality $
\mathrm{Tr\,}f(A+B) \le \mathrm{Tr\,}f(A)+f(B) 
$
 combined with $f(0)=0$ means that $f(t)$ is concave. This shows that (1.2) is more subtle than (1.1).
\end{remark}

\vskip 5pt
\begin{remark} There exists a norm version of Rotfel'd inequality which considerably improves (3.1). {\it If $f:[0,\infty)\to[0,\infty)$ is concave and $A, B\in\bM_n^+$, then
$$
\| f(A+B) \| \le  \|f(A) +f(B) \|
$$
for all symmetric norms.} The case of operator concave functions is given in \cite{AZ} and the general case is established in \cite{BU}, see also \cite{BL} for further results. Concerning differences, the following inequality holds
$$
\| f(A)-f(B) \| \le \| f(|A-B|) \|
$$
for all symmetric norms,  $A,B\in\bM_n^+$, and  operator monotone  functions $f:[0,\infty)\to[0,\infty)$. This is a famous result of Ando \cite{Ando}. Here the operator monotony assumption is essential, see \cite{AuAu} for some counterexamples. A very interesting paper by Mathias \cite{M} gives a direct proof, without using the integral representation of operator monotone functions.
\end{remark}

\vskip 5pt
\begin{remark} There exists also some subaditivity results involving  convex functions \cite{BH1}. For instance:
{\it
Let $g(t)=\sum_{k=0}^m a_kt^k$ be a polynomial of degree $m$ with all non-negative
coefficients. Then, for all positive operators $A,\,B$ and all symmetric norms,}
\begin{equation*}
\| g(A+B) \|^{1/m} \le \|g(A) \|^{1/m} + \| g(B) \|^{1/m}.
\end{equation*} 
\end{remark}

\vskip 5pt
\begin{remark} It is not known wether the monotony assumption in Theorem 3.1 can be deleted.
\end{remark}


\vskip 20pt
J.-C. Bourin,

Laboratoire de math\'ematiques,

Universit\'e de Franche-Comt\'e,

25 000 Besancon, France.

jcbourin@univ-fcomte.fr

\vskip 20pt
Eun-Young Lee

 Department of mathematics,

Kyungpook National University,

 Daegu 702-701, Korea.

eylee89@ knu.ac.kr

\vskip 20pt\noindent

\end{document}